\documentclass[reqno]{amsart}
\usepackage{amsfonts}
\usepackage{amssymb}
\usepackage{latexsym}
\usepackage{amssymb}
\usepackage{amsmath}
\usepackage{mathrsfs}
\usepackage{amsbsy}
\usepackage{cite}
\usepackage{hyperref}
\usepackage[latin1]{inputenc}

\newtheorem{theorem}{Theorem}

\newtheorem{proposition}{Proposition}

\begin{document}
\title[Class-preserving automorphisms of universal hyperlinear groups]{%
Class-preserving automorphisms\\
of universal hyperlinear groups}
\author{Martino Lupini}
\address{Department of Mathematics and Statistics\\
N520 Ross, 4700 Keele Street, M3J 1P3\\
Toronto (Ontario), Canada}
\email{mlupini@mathstat.yorku.ca}
\thanks{The author was supported by the York University Elia Scholars Program%
}
\date{}
\dedicatory{}

\begin{abstract}
We show that the group of class-preserving automorphisms of a universal
hyperlinear group has index $2$ inside the group of all automorphisms.
\end{abstract}

\maketitle

\section{Introduction\label{Section:introduction}}

Suppose that $G$ is a group. An automorphism $\alpha $ of $G$ is \emph{%
class-preserving} if for every $g\in G$ the elements $g$ and $\alpha (g)$ of 
$G$ belong to the same conjugacy class. It is easily observed that
class-preserving automorphisms of $G$ form a normal subgroup of the group of
automorphisms of $G$. It is shown in \cite{Paunescu} that if $\mathcal{U}$
is a nonprincipal ultrafilter over $\mathbb{N}$, then every automorphism of
the ultraproduct $\prod\nolimits_{\mathcal{U}}S_{n}$ of the finite symmetric
groups endowed with the normalized Hamming metric (as defined in \cite[%
Section 2.4]{pestov}) is class-preserving. The groups $\prod\nolimits_{%
\mathcal{U}}S_{n}$ are sometimes called \emph{universal sofic groups}, since
they contain any countable discrete sofic group as a subgroup, see \cite[%
Section 3]{pestov}.

In this note we consider \emph{universal hyperlinear groups}, i.e.\
ultraproducts $\prod\nolimits_{\mathcal{U}}U_{n}$ of the finite rank unitary
groups endowed with the normalized Hilbert-Schmidt metric, as defined in 
\cite[Section 2.4]{pestov}. We show that for any nonprincipal ultrafilter $%
\mathcal{U}$ over $\mathbb{N}$ the group of class-preserving automorphism of 
$\prod\nolimits_{\mathcal{U}}U_{n}$ has index $2$ in the automorphism group
of $\prod\nolimits_{\mathcal{U}}U_{n}$. The same statement holds when one
replaces $\prod\nolimits_{\mathcal{U}}U_{n}$ with the unitary group $U(%
\mathcal{R}^{\mathcal{U}})$ of the ultrapower $\mathcal{R}^{\mathcal{U}}$ of
the separable hyperfinite II$_{1}$ factor (an introduction to II$_{1}$
factors and a definition of the hyperfinite II$_{1}$ factor can be found in 
\cite[Section III.1]{Blackadar}). It is worth observing that by \cite[%
Proposition 2.4.6]{Capraro-Lupini} if the Continuum Hypothesis holds, then $%
U(\mathcal{R}^{\mathcal{U}})$ and $\prod\nolimits_{\mathcal{U}}U_{n}$ have
outer automorphisms (and in fact $2^{\aleph _{1}}$ many of them).


The rest of this note is divided into two sections: In Section \ref{Section:
locally inner} we recall a fact about automorphisms of countably saturated II%
$_{1}$ factors whose proof is essentially contained in \cite{Sherman-notes};
In Section \ref{Section: pointwise inner} we present a proof of the main
theorem, based on Theorem 2 from \cite{Dye2}. In the following $\mathcal{U}$
will always be assumed to be a fixed nonprincipal ultrafilter over $\mathbb{N%
}$. The separable hyperfinite II$_{1}$ factor will be denoted as customary
by $\mathcal{R}$. The set of natural numbers $\mathbb{N}$ will be assumed
not to contain $0$, and a positive real number will be assumed to be
strictly greater than zero. A natural number $n$ will be regarded as a
finite ordinal and hence identified with the set $\left\{ 0,1,\ldots
,n-1\right\} $ of its predecessors.

\section{Automorphisms of countably saturated \texorpdfstring{II$_{1}$}{II1}
factors\label{Section: locally inner}}

In this section we make use of terminology and results from the logic for
metric structures. (An introduction to this subject can be found in \cite%
{BBHU}.) In particular we consider II$_{1}$ factors as structures in the
language of tracial von Neumann algebras, as described in \cite[Section 2.3.2%
]{FHS2}.

Suppose that $M$ is a II$_{1}$ factor, $\beta $ is an automorphism of $%
\mathcal{M}$, and $a$ is a normal element in the unit ball of $\mathcal{M}$.
Approximating $a$ by normal elements with finite spectrum, it is easy to see
that there is a sequence $\left( u_{n}\right) _{n\in \mathbb{N}}$ of unitary
elements of $\mathcal{M}$ such that%
\begin{equation*}
\lim_{n\rightarrow +\infty }\left\Vert \beta (a)-u_{n}au_{n}^{\ast
}\right\Vert _{2}=0
\end{equation*}%
where $\left\Vert \cdot \right\Vert _{2}$ is the Hilbert-Schmidt norm
associated with the unique trace of $\mathcal{M}$. Thus the formula (with
parameters)%
\begin{equation*}
\psi (z)\equiv \left\Vert \beta (a)-zaz^{\ast }\right\Vert _{2}+\left\Vert
zz^{\ast }-1\right\Vert _{2}+\left\Vert z^{\ast }z-1\right\Vert _{2}
\end{equation*}%
is approximately realized in $\mathcal{M}$, i.e.\ for every positive real
number $\varepsilon $ there is an element $u$ in the unit ball of $\mathcal{M%
}$ such that $\psi (u)<\varepsilon $. If moreover $\mathcal{M}$ is countably
saturated as in \cite[Section 4.4]{FHS2}, then the formula $\psi $ is
actually realized in $\mathcal{M}$, i.e.\ there is a (necessarily unitary)
element $u$ of $\mathcal{M}$ such that $\psi (u)=0$ and hence%
\begin{equation*}
uau^{\ast }=\beta (a)\text{.}
\end{equation*}%
This concludes the proof of Proposition \ref{Proposition: locally inner}.

\begin{proposition}
\label{Proposition: locally inner}Suppose that $\mathcal{M}$ is a countably
saturated {\upshape II}$_{1}$ factor, and $\beta $ is an automorphism of $%
\mathcal{M}$. If $a$ is a normal element of $\mathcal{M}$, then $a$ and $%
\beta (a)$ are conjugate by a unitary element of $\mathcal{M}$.
\end{proposition}

Recall that by \cite[Proposition 4.11]{FHS2} an ultraproduct of a sequence
of tracial von Neumann algebras with respect to a nonprincipal ultrafilter
is countably saturated.

Proposition \ref{Proposition: locally inner} can also be proved using
methods and results from \cite{Sherman-notes}. In fact it is not difficult
to see that the same proof as in \cite{Sherman-notes} shows that Theorem 3.1
and Corollary 3.6 from \cite{Sherman-notes} hold not only for ultrapowers of
II$_{1}$ factors, but also more generally for any countably saturated II$%
_{1} $ factor. If now $a$ is a normal element of $\mathcal{M}$, and $A$ is a
the C*-subalgebra of $\mathcal{M}$ generated by $a$, then $A$ is abelian
and, in particular, nuclear. It therefore follows from the described
generalization of \cite[Corollary 3.6]{Sherman-notes} that the restriction
of $\beta $ to $A $ coincides with the restriction of some inner
automorphism of $\mathcal{M} $. In particular $a$ and $\beta (a)$ are
conjugate by a unitary element of $\mathcal{M}$.


\section{Automorphisms of \texorpdfstring{$U(\mathcal{R}^{%
\mathcal{U}})$}{U(RU)}}

\label{Section: pointwise inner} The main result of this section is Theorem %
\ref{Theorem: pointwise inner}. Recall that an automorphism $\alpha $ of a
group $G$ is class-preserving if for every $x\in G$ the elements $x$ and $%
\alpha (x)$ of $G$ are conjugate. The class-preserving automorphisms of $G$
form a normal subgroup of the automorphisms group of $G$.

\begin{theorem}
\label{Theorem: pointwise inner}The group of class-preserving automorphisms
of $U(\mathcal{R}^{\mathcal{U}})$ has index $2$ in the automorphisms group
of $U(\mathcal{R}^{\mathcal{U}})$. Moreover if $\alpha $ is any automorphism
of $U(\mathcal{R}^{\mathcal{U}})$ then $\alpha $ is an isometry with respect
to the distance induced by the Hilbert-Schmidt norm on $\mathcal{R}^{%
\mathcal{U}}$. The same statement holds for automorphisms of the
ultraproduct $\prod\nolimits_{\mathcal{U}}U_{n}$ of the sequence of unitary
groups endowed with the normalized Hilbert-Schmidt metric.
\end{theorem}

Identify $\mathcal{R}$ with the von Neumann algebra tensor product of
infinitely many copies of $\mathbb{M}_{2}$. If $a\in \mathbb{M}_{n}$ define $%
\bar{a}$ the element of $\mathbb{M}_{2}$ obtained replacing every entry of $%
a $ with the corresponding complex conjugate. It is immediate to verify that
the function $a\mapsto \bar{a}$ is a conjugate linear automorphism of $%
\mathbb{M}_{2}$. Moreover the function%
\begin{equation*}
a_{0}\otimes \cdots \otimes a_{n-1}\mapsto \bar{a}_{0}\cdots \otimes \bar{a}%
_{n-1}
\end{equation*}%
induces a conjugate linear automorphism of $\mathcal{R}$. Passing to the
ultrapower one obtains a conjugate linear automorphism $\gamma $ of $%
\mathcal{R}^{\mathcal{U}}$ preserving the Hilbert-Schmidt norm. The
restriction of $\gamma $ to $U(\mathcal{R}^{\mathcal{U}})$ defines an
automorphism of $U(\mathcal{R}^{\mathcal{U}})$ that is not class-preserving
but it is an isometry with respect to the distance induced by the
Hilbert-Schmidt norm on $\mathcal{R}^{\mathcal{U}}$. We will show in the
following that if $\alpha $ is any automorphism of $U(\mathcal{R}^{\mathcal{U%
}})$, then either $\alpha $ or $\gamma \circ \alpha $ is class-preserving.
This in particular will show that $\alpha $ is an isometry with respect to
the distance induced by the Hilbert-Schmidt norm on $\mathcal{R}^{\mathcal{U}%
}$.

Suppose that $\alpha $ is an automorphism of $U(\mathcal{R}^{\mathcal{U}})$.
Recall that by cornerstone results of Dixmier-Lance \cite{Dix-Lance} and
McDuff \cite{McD} an ultrapower of a II$_{1}$ factor is a II$_{1}$ factor
(this also follows from the fact that the class of II$_{1}$ factors is
axiomatizable, see \cite[Proposition 3.4]{BBHU}). Thus $\mathcal{R}^{%
\mathcal{U}}$ is a II$_{1}$ factor, which is moreover countably saturated by 
\cite[Proposition 4.11]{FHS2}. By Theorem 2 from \cite{Dye2} together with
the result of Broise from \cite{Broise} that the unitary group of a II$_{1}$
factor does not have characters, there is a linear or conjugate linear
*-isomorphism $\beta :\mathcal{R}^{\mathcal{U}}\rightarrow \mathcal{R}^{%
\mathcal{U}}$ whose restriction to $U(\mathcal{R}^{\mathcal{U}})$ is $\alpha 
$. If $\beta $ is a linear *-isomorphism then by Proposition \ref%
{Proposition: locally inner} for every $u\in U(\mathcal{R}^{\mathcal{U}})$
there is $v\in U(\mathcal{R}^{\mathcal{U}})$ such that%
\begin{equation*}
vuv^{\ast }=\beta (u)=\alpha (u)\text{.}
\end{equation*}%
and hence $\alpha $ is class-preserving. If $\beta $ is a conjugate linear
*-isomorphism, then $\gamma \circ \beta $ is a linear *-isomorphism whose
restriction to $U(\mathcal{R}^{\mathcal{U}})$ is $\gamma \circ \alpha $.
Reasoning as before one concludes that $\gamma \circ \alpha $ is
class-preserving. This concludes the proof Theorem \ref{Theorem: pointwise
inner}.

In order to prove the analogous statement for $\prod\nolimits_{\mathcal{U}%
}U_{n}$ observe that $\prod\nolimits_{\mathcal{U}}U_{n}$ can be identified
with the unitary group of $\prod\nolimits_{\mathcal{U}}\mathbb{M}_{n}$ by 
\cite[Proposition 2.1]{Ge-Hadwin}, see also \cite[Exercise II.9.6]%
{Capraro-Lupini}. One can thus run the same argument where $\mathcal{R}^{%
\mathcal{U}}$ is replaced by $\prod\nolimits_{\mathcal{U}}\mathbb{M}_{n}$
which is again a countably saturated II$_{1}$ factor by \cite[Proposition
4.11]{FHS2}.

\subsection*{Acknowledgments}
We would like to thank Goulnara Arzhantseva for her comments, and for pointing out a mistake in the first version of these notes; We are also grateful to Valerio Capraro, \L ukasz Grabowski, and Liviu P\u aunescu for their remarks and suggestions.















\end{document}